\documentclass[leqno]{amsart}
\usepackage{amssymb, amsmath, amsthm}
\usepackage[all]{xy}
\usepackage{color}

\newtheorem{thm}{Theorem}[section]
\newtheorem{theorem}[thm]{Theorem}

\newtheorem{prop}[thm]{Proposition}

\newtheorem{lemma}[thm]{Lemma}

\theoremstyle{definition} 
\newtheorem{defn}[thm]{Definition}

\newcommand{\cee}{\ensuremath{\mathbb C}}

\newcommand{\Mnd}{{\sf Mnd} }

\newcommand{\tensor}{\otimes}

\newcommand{\teu}{(T, \eta, \mu)}
\newcommand{\C}{\mathbb {C}}

\newcommand{\sdelta}{\ensuremath{2  \Delta}}

\begin{document}
\date{\today}
\title[2 Monads and Simplicial Structures]{ General 2-Dimensional Adjunctions, Universal Monads and Simplicial Structures }
\author{John Lauchlin MacDonald and Laura Scull}
\begin{abstract} We use the general notion of 2-dimensional adjunction with given coherence equations as introduced by MacDonald-Stone \cite{8}, building on earlier work by  Gray \cite{4}, to derive coherence equations for a general 2-monad, which we refer to as a lax-Gray monad in keeping with terminology from \cite{2, 7}.  The free lax-Gray 2-monad on one object,  $\sf 2-Mnd$,  may be regarded as the suspension of a lax 2-dimensional analogue of the simplicial category $\Delta$.  We   call this analogue $\Delta_{LG}$ for lax Gray $\Delta$.  This is analogous to the way that the free 1-monad {\sf Mnd} (as presented in Schanuel-Street \cite{10}) is a concrete example of the suspension of the simplicial category $\Delta$ described by MacLane \cite{9}.    \end{abstract}
\maketitle

\section{Introduction}  In this paper, we generalize the classical simplicial category $\Delta$ to a higher dimensional analogue object.  To do this, we interpret the category $\Delta$ as the arrow category of a 2-category ${\sf 1-Mnd}$ created by the universal monad.  This category has   one object, where the arrow category is freely generated by one 1-cell $T$ and two 2-cells,  $\eta:  1 \to T$ and $\mu:  T^2 \to T$.  These generate the classic face and degeneracy maps of $\Delta$ and satisfy the standard coherence  equations $\mu . T\mu = \mu . \mu T$ and $\mu.T\eta = 1 = \mu.\eta T$.   Thus the category $\Delta$ can be seen as a desuspension of the  universal monad category.  

We will generalize this by creating a universal 2-monad category $\sf {2-Mnd}$, and then defining our higher dimensional simplicial category $\Delta_{LG}$ as a desuspension of  this category.   To create the universal 2-monad, we will 
 use the general notion of 2-dimensional adjunction with given coherence equations as introduced by MacDonald-Stone \cite{8}, building on earlier work by  Gray \cite{4}.  We will work in the setting of lax-Gray categories, and  derive coherence equations for a general 2-monad, which we refer to as a lax-Gray monad in keeping with terminology from \cite{2, 7}.  The free lax-Gray 2-monad on one object,  $\sf 2-Mnd$,  will give us the  suspension of a lax 2-dimensional analogue of the simplicial category $\Delta$.  We note that the work by Lack \cite{7} on pseudomonads contains a related suspension pseudomonad. 

\section{Background:  Adjunctions and Coherence} 

In this section we discuss how adjunctions are generalized to higher dimensions, and set the context which we will work in.  
To set notation, we will picture an isomorphism with the following diagram \[ \xymatrix{ (1) & X = \ar@<+1ex>[r]^f &  = A\ar@<1ex>[l]^u } \] 
indicating that the composites $fu$ and $uf$ are the identities on $A$ and $X$ respectively.  Similarly, an adjunction is denoted with the diagram 
 \[ \xymatrix{ (2) & X \overset{\eta}{\Rightarrow} \ar@<+1ex>[r]^F&  \overset{\epsilon}{\Rightarrow} A\ar@<1ex>[l]^U } \] where the equations $1_X = uf$ and $fu = 1_A$ of Equation (1) are replaced by 2-cells $\eta:  1_X \Rightarrow UF$ and $\epsilon:  FU \Rightarrow 1_A$.  These 2-cells must satisfy their own coherence equations on the next level,  $ \epsilon F. F \eta = 1_F$  and $1_U = U \epsilon.\eta U$.

Taking this one step further,  we consider here what happens when the coherence conditions of the adjunction in (2), given by the equations $\epsilon F . F \eta = 1_F$ and $1_U = U\epsilon . \eta U$,   are themselves replaced by 3-cells with their own coherence conditions.    In this case, as in MacDonald-Stone \cite{4}, we have a general 2-dimensional adjunction pictured by 
  \[ \xymatrix{ (3) &  X \overset{\eta}{\Rightarrow} \ar@<+1ex>[r]^{f \Downarrow \, \,  F}&  \overset{\epsilon}{\Rightarrow} A\ar@<1ex>[l]^{ u \Downarrow \, \, U} } \]  
extending diagram (2).  Explicitly, this is defined by:

\begin{defn} A 2-dimensional adjunction is defined by a  6-tuple $(F, U, \eta, \epsilon, f, \mu)$, where
  \begin{enumerate}
  \item[]  $F:  X \to A $ and $U:  A \to X$ are strict 2-functors
  \item[] $\eta:  X \Rightarrow UF$ and $\epsilon:  FU \Rightarrow A$ are lax natural transformations 
  \item[] $f:  \epsilon F. F \eta  \Rrightarrow 1_F $ and $u: 1_U  \Rrightarrow U\eta. \epsilon U$ are modifications between lax transformations. 
    \end{enumerate}
    and $f, u$  satisfy the following coherence conditions
    
    \begin{enumerate}
    \item[{[$1_\eta$]}]  $Uf . \eta .. U \epsilon F .  \eta \eta  ..  \mu F .  \eta = 1 _\eta$
    
      \item[{[$1_\epsilon$]}]  $ 1_\epsilon =\epsilon .  fU  ..  \epsilon \epsilon .  F \eta U  ..  \epsilon .  F\mu $
    \end{enumerate} where $..$ denotes composition of modifications.  
  \end{defn}
  
 Each of these equations can be presented geometrically as an equality between composition of modifications on the front and back sides of a cube

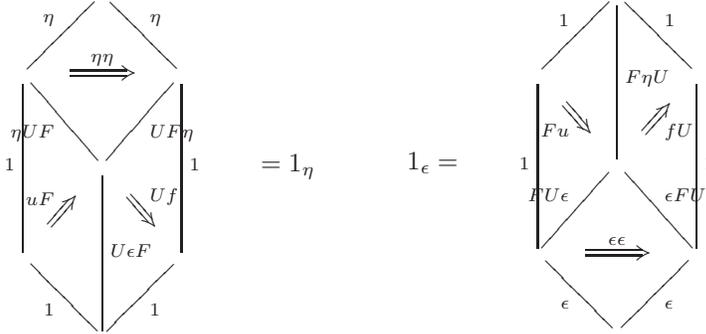
\[ \xymatrix{ &  \ar@{-}[dl]_{\eta} \ar@{-}[dr]^{\eta} & \\
 \ar@{-}[dr]_{\eta UF}  \ar@{}[rr]^{\eta \eta}^(.25){}="a"^(.75){}="b" \ar@{=>} "a";"b"  \ar@{-}[dd]_1 &&  \ar@{-}[dd]^1 \ar@{-}[dl]^{UF \eta} \\
&   \ar@{-}[dd]^{U \epsilon F}   \ar@{}[dr]^{Uf}^(.25){}="c"^(.75){}="d" \ar@{=>} "c";"d"  &  & =  1_\eta\\
 \ar@{}[ur]^{uF}^(.25){}="g"^(.75){}="f" \ar@{=>} "g";"f"  \ar@{-}[dr]_1 &&   \ar@{-}[dl]^1  \\
&  & \\
} \, \, \phantom{www} \, \,\xymatrix{&  &  \ar@{-}[dl]_{1} \ar@{-}[dr]^{1}  \ar@{-}[dd]^{F \eta U }& \\
&    \ar@{}[dr]_{Fu}^(.25){}="c"^(.75){}="d" \ar@{=>} "c";"d"  \ar@{-}[dd]_1    &&  \ar@{-}[dd]^1  \\
  1_\epsilon = & &    \ar@{}[ur]_{fU}^(.25){}="g"^(.75){}="f" \ar@{=>} "g";"f"    \ar@{-}[dl]_{FU \epsilon} \ar@{-}[dr]^{\epsilon FU }&  & \\
&  \ar@{-}[dr]_\epsilon \ar@{}[rr]^{\epsilon \epsilon}^(.25){}="a"^(.75){}="b" \ar@{=>} "a";"b"  &&   \ar@{-}[dl]^\epsilon \\
&&   & \\
}  \]

Now if we compare these various definitions, we see that the definition (1) of isomorphism makes sense in any category, including particularly the prototype category of sets and functions.  Similarly diagram (2) defining an adjunction makes sense in any 2-category as in Kelly-Street \cite{10}.  Typically adjunctions are first presented in the prototype 2-category {\sf Cat} of categories, functors and natural transformations. 

It might seem natural, then, to generalize to 3-categories to consider definition (3).  However, in view of earlier work by MacDonald-Stone \cite{4} and other current literature \cite{1, 2, 7}, it seems that actually the right direction is to move to Gray-categories and lax-Gray categories, which provide the flexibility needed for diagram (3) and its coherence equations.  

According to nLab:  {\it The prototypical Gray-category is {\sf Gray}, which consists of strict 2-categories, strict 2-functors, pseudonatural transformations and modifications.}  In this paper, we will work in lax-Gray categories and particularly in its prototypical category {\sf lax-Gray} consisting of strict 2-categories, strict 2-functors, lax natural transformations and modifications.

{\sf Gray} has a tensor product $B \otimes C$ defined by ${\sf 2-Cat}(B \otimes C, D) \equiv {\sf 2-Cat}(B, P_S(C, D))$ where $P_S(C, D)$ denotes the 2-category of 2-functors, pseudonatural transformations and modifications from $C$ to $ D$.  Similarly, {\sf lax-Gray} has a tensor product $B \otimes_L C$ defined by ${\sf 2-Cat}(B \otimes_L C, D) \equiv {\sf 2-Cat}(B, P_L(C, D))$ where $P_L(C, D)$ denotes the 2-category of 2-functors, lax natural transformations and modifications from $C$ to $ D$.  In this sense, we can recast Gray-categories as categories enriched over {\sf Gray}, that is, $V$-categories for $V = {\sf Gray}$.  Similarly,  lax-Gray categories are $V$-categories for V = {\sf lax-Gray}. 

This is the framework in which we will interpret the general 2-dimensional adjunction of (3) described by equations $[1_\eta] $ and $[1_\epsilon]$.   These structures were described earlier in MacDonald-Stone \cite{4} where they were presented geometrically and called soft adjunctions.  Here we call them lax-Gray adjunctions to bring our terminology in line with recent work \cite{1, 2, 5, 7}. 

To set notation, the 2-cell direction that we use for a lax natural transformation $\varphi:  S \to T:  \mathbb{A} \to \mathbb{B}$ is given by 
 \[ \xymatrix{ SA \ar[r]^{\varphi_A} \ar[d]_{Sa} & TA \ar[d]^{Ta}  \\ 
      SB \ar[r]_{\varphi_B} \ar@{}[ur]^{\varphi_a}^(.3){}="a"^(.7){}="b"  \ar@{=>} "a";"b"  & TB } \] 
      for a morphism $a:  A \to B$ in $\mathbb{A}$ with usual equations $\varphi_{1_\mathbb{A}} = 1 _{\varphi_A}, \varphi_{ba} = T_b\varphi_a . \varphi_b S_a$ and $T \zeta . \varphi_a .. \varphi_a = \varphi_{a'} .. \varphi_B . S\zeta$ for a 2-cell $\zeta:  a \Rightarrow a'$.  In particular,  the lax natural transformation $\epsilon FU \to 1:  \mathbb{A} \to \mathbb{A}$ of  (3) is defined by:  for each 1-cell a of $\mathbb{A}$, there is  a 2-cell
 \[ \xymatrix{ FUA \ar[r]^{\epsilon_A} \ar[d]_{FUa} & A \ar[d]^{a}  \\ 
      FUB \ar[r]_{\epsilon_B} \ar@{}[ur]^{\epsilon_a}^(.3){}="a"^(.7){}="b"  \ar@{=>} "a";"b"  & B } \] 
and consequently for each object $A$ of $\mathbb{A}$, there is a 2-cell 
 \[ \xymatrix{ FUFUA \ar[r]^{\epsilon_{FUA}} \ar[d]_{FU\epsilon_A} & FUA \ar[d]^{\epsilon_A}  \\ FUA \ar[r]_{\epsilon_A} \ar@{}[ur]^{\epsilon_{\epsilon A}}^(.3){}="a"^(.7){}="b"  \ar@{=>} "a";"b"  & A } \] 
which determines the $\epsilon \epsilon $ of the coherence condition $[1_\epsilon]$ for diagram (3).   In particular, if we let $A = FX$ then we get  
 \[ \xymatrix{ FUFUFX \ar[rr]^{\epsilon_{FUFX}} \ar[d]_{FU\epsilon_{FX}} &&  FUFX \ar[d]^{\epsilon_{FX}}  \\ FUFX \ar[rr]_{ \epsilon_{FX}} \ar@{}[urr]^{\epsilon{\epsilon _{FX}}}^(.35){}="a"^(.65){}="b"  \ar@{=>} "a";"b"  &&  FX } \] 

If we apply U to this diagram and let $T =UF, \mu = U\epsilon F$ and $a= U\epsilon \epsilon F$ then we have the 2-cell
\[ \xymatrix{ T^3X \ar[rr]^{\mu_{TX}} \ar[d]_{T\mu_{X}} &&  T^2X \ar[d]^{\mu_X}  \\ T^2X \ar[rr]_{ \mu_{X}} \ar@{}[urr]^{a_X}^(.35){}="a"^(.65){}="b"  \ar@{=>} "a";"b"  &&  TX } \] 
This cell $a_X$ is an example of one used in the coherence conditions for a general 2-monad.  

It is shown in MacDonald-Stone \cite{8} that this equational, or coherence, form for the definition of a lax-Gray adjunction is equivalent to a description by adjunction of hom-categories as well as to one by universal properties, in analogy with the 1-dimensional case described by MacLane \cite{9}.

\section{Lax-Gray Monads}

A general 2-dimensional monad, or lax-Gray monad, on a 2-category $X$ is a 6-tuple $ (T, \eta, \mu, a, p, q) $ where $\eta:  1 \to T$ and $\mu:  T^2 \to T$ are lax natural transformations, and $a:  \mu . T\mu \Rightarrow \mu . \mu T, \, \, \,  p:  \mu . T \eta \to 1_T$ and $q:  1_T \to \mu . \eta T$ are modifications.  This data may be partially displayed in the diagram 
 \[ \xymatrix{   1 \ar[r]^{\eta} & \, \,  T \, \, \, \, \ar@<+1ex>@{=>}[r]^q &  \ar@<+1ex>@{=>}[l]^p  \ar@<+3ex>[r]^{\eta T} \ar@<-3ex>[r]^{{T\eta}}&  \, \,  T^2  \, \,  \ar@<0ex>[l]_\mu \ar@<+4ex>@{=>}[r]^{qT}  \ar@<-2ex>@{=>}[r]^{Tq} &  \ar@<-1ex>@{=>}[l]_{pT} \ar@<+5ex>@{=>}[l]_{Tp}  \ar@<+5ex>[r]^{\eta T^2}    \ar@<0ex>[r]^{T\eta T}    \ar@<-5ex>[r]^{ T^2 \eta}&   \, \,  T^3   \, \,   \ar@<+2.5ex>[l]_{T\mu }  \ar@<-2.5ex>[l]_{ \mu T}  \\  } \]  
subject to the following five coherence axioms.

[M$a$]   \[a . \mu T^2  ..  \mu .  \mu \mu  ..  a .  T^2 \mu = \mu  ..   aT  ..  a .  T\mu T  ..  \mu .  Ta\] 
We can present this as a commuting cubic hexagon of modifications where the composition of the front three sides equals the composition on the back three sides.  We abbreviate this as 
\[ \xymatrix{ & \ar@{-}[dr]^{\mu T^2} \ar@{-}[dl]_{T^2 \mu} && & & \ar@{-}[dd] \ar@{-}[dr]^{\mu T^2 }  \ar@{-}[dl]_{ T^2 \mu} & &\\ \ar@{-}[dd]_{T\mu} \ar@{-}[dr]_{\mu T}  \ar@{}[rr]^{\mu \mu }^(.3){}="a"^(.7){}="b"  \ar@{=>} "a";"b"  &  &\ar@{-}[dd]^{\mu T}  \ar@{-}[dl]^{T\mu} &&  \ar@{-}[dd]_{T \mu}  \ar@{}[dr]_{Ta}^(.25){}="l"^(.75){}="m"   \ar@{=>} "l";"m"&  \ar@{=>} "l";"m" & \ar@{-}[dd]^{\mu T}   &   &  \\
& \ar@{-}[dd]  \ar@{}[dr]^{a}^(.25){}="f"^(.77){}="g"  \ar@{=>} "f";"g" & & = &  &  \ar@{-}[dr]  \ar@{-}[dl]  \ar@{}[ur]_{aT}^(.25){}="s"^(.75){}="r"   \ar@{=>} "s";"r" & &  \\ 
 \ar@{-}[dr]_{\mu}  \ar@{}[ur]^{a}^(.25){}="c"^(.75){}="d"  \ar@{=>} "c";"d" && \ar@{-}[dl]^{\mu}   &&   \ar@{-}[dr]_{\mu} \ar@{}[rr]^{a }^(.3){}="j"^(.7){}="k"  \ar@{=>} "j";"k"&&  \ar@{-}[dl]^{\mu}  \\ &&&&&&&& \\ 
} \]  

[M$\eta$]  \[ p .  \eta  ..  \mu .  \eta \eta  ..  q . \eta =    1_{\eta}\]

\[ \xymatrix{ & \ar@{-}[dr]^{\eta} \ar@{-}[dl]_{\eta} && & & \ar@{-}[dd] \ar@{-}[dr]^{\eta}  \ar@{-}[dl]_{\eta} & &\\ \ar@{-}[dd]_{1} \ar@{-}[dr]  \ar@{}[rr]^{\eta \eta  }^(.3){}="a"^(.7){}="b"  \ar@{=>} "a";"b"  &  &\ar@{-}[dd]^{1}  \ar@{-}[dl] &&  \ar@{-}[dd]_{1}  \ar@{}[dr]_{1}^(.25){}="l"^(.75){}="m"   \ar@{=>} "l";"m"&  \ar@{=>} "l";"m" & \ar@{-}[dd]^{1}   &   &  \\
& \ar@{-}[dd]  \ar@{}[dr]^{p}^(.25){}="f"^(.77){}="g"  \ar@{=>} "f";"g" & & = &  &  \ar@{-}[dr]  \ar@{-}[dl]  \ar@{}[ur]_{1}^(.25){}="s"^(.75){}="r"   \ar@{=>} "s";"r" & &  \\ 
 \ar@{-}[dr]^{1}  \ar@{}[ur]^{q}^(.25){}="c"^(.75){}="d"  \ar@{=>} "c";"d" && \ar@{-}[dl]_{1}  &&   \ar@{-}[dr]_{1} \ar@{}[rr]^{1}^(.3){}="j"^(.7){}="k"  \ar@{=>} "j";"k"&&  \ar@{-}[dl]^{1}  \\ &&&&&&&& \\ 
} \]

[M$\mu$] \[ 1 _{\mu} = \mu  .  pT   ..  a .  T \eta T  ..  \mu . Tq \]
\[ \xymatrix{ & \ar@{-}[dr]^{1} \ar@{-}[dl]_{1} && & & \ar@{-}[dd] \ar@{-}[dr]^{1}  \ar@{-}[dl]_{1} & &\\ \ar@{-}[dd]_{1} \ar@{-}[dr]  \ar@{}[rr]^{1 }^(.3){}="a"^(.7){}="b"  \ar@{=>} "a";"b"  &  &\ar@{-}[dd]^{1}  \ar@{-}[dl] &&  \ar@{-}[dd]_{1}  \ar@{}[dr]_{Tq}^(.25){}="l"^(.75){}="m"   \ar@{=>} "l";"m"&  \ar@{=>} "l";"m" & \ar@{-}[dd]^{1}   &   &  \\
& \ar@{-}[dd]  \ar@{}[dr]^{1}^(.25){}="f"^(.77){}="g"  \ar@{=>} "f";"g" & & = &  &  \ar@{-}[dr]  \ar@{-}[dl]  \ar@{}[ur]_{pT}^(.25){}="s"^(.75){}="r"   \ar@{=>} "s";"r" & &  \\ 
 \ar@{-}[dr]_{\mu}  \ar@{}[ur]^{1}^(.25){}="c"^(.75){}="d"  \ar@{=>} "c";"d" && \ar@{-}[dl]^{\mu}   &&   \ar@{-}[dr]_{\mu} \ar@{}[rr]^{a }^(.3){}="j"^(.7){}="k"  \ar@{=>} "j";"k"&&  \ar@{-}[dl]^{\mu}  \\ &&&&&&&& \\ 
} \]  

[M$q$]   \[ a .  \eta T^2  ..  \mu . \eta \mu  ..  q.  \mu = \mu .  qT \]
\[ \xymatrix{ & T^2 \ar@{-}[dr]^{\eta T^2 } \ar@{-}[dl]_{\mu} && & & T^2  \ar@{-}[dd]_{1} \ar@{-}[dr]^{\eta T^2}  \ar@{-}[dl]_{\mu } & &\\  T \ar@{-}[dd]_{1} \ar@{-}[dr]_{\eta T}  \ar@{}[rr]^{\eta  \mu }^(.3){}="a"^(.7){}="b"  \ar@{=>} "a";"b"  &  & T^3 \ar@{-}[dd]^{\mu T}  \ar@{-}[dl]^{T \mu}  && T \ar@{-}[dd]_{1}  \ar@{}[dr]_{1}^(.25){}="l"^(.75){}="m"   \ar@{=>} "l";"m"&  \ar@{=>} "l";"m" & \ar@{-}[dd]^{\mu T}    T^3 & &  \\
& \ar@{-}[dd]^{\mu}   \ar@{}[dr]^{a}^(.25){}="f"^(.77){}="g"  \ar@{=>} "f";"g" & & = &  & T^2 \ar@{-}[dr]^{1}   \ar@{-}[dl]_{\mu}  \ar@{}[ur]_{aT}^(.25){}="s"^(.75){}="r"   \ar@{=>} "s";"r" & &  \\ 
T  \ar@{-}[dr]_{1}  \ar@{}[ur]^{q}^(.25){}="c"^(.75){}="d"  \ar@{=>} "c";"d" && T^2 \ar@{-}[dl]^{\mu}   &&  T  \ar@{-}[dr]_{1} \ar@{}[rr]^{1 }^(.3){}="j"^(.7){}="k"  \ar@{=>} "j";"k"&& T^2  \ar@{-}[dl]^{\mu}  \\ & T &&&&T && \\ 
} \]

[M$p$]   \[ p .  \mu  ..  T\mu .  \mu \eta  ..  a .  T^2 \eta= \mu .  Tp \]
\[ \xymatrix{ & T^2 \ar@{-}[dr]^{\mu } \ar@{-}[dl]_{T^2 \eta} && & & T^2  \ar@{-}[dd]_{1} \ar@{-}[dr]^{\mu}  \ar@{-}[dl]_{T^2 \eta } & &\\  T^3 \ar@{-}[dd]_{T \mu} \ar@{-}[dr]_{\mu T}  \ar@{}[rr]^{\mu \eta   }^(.3){}="a"^(.7){}="b"  \ar@{=>} "a";"b"  &  & T \ar@{-}[dd]^{1}  \ar@{-}[dl]^{T \eta}  && T^3 \ar@{-}[dd]_{T \mu}  \ar@{}[dr]_{Tp}^(.25){}="l"^(.75){}="m"   \ar@{=>} "l";"m"&  \ar@{=>} "l";"m" & \ar@{-}[dd]^{1}    T & &  \\
& \ar@{-}[dd]^{T \mu }   \ar@{}[dr]^{p}^(.25){}="f"^(.77){}="g"  \ar@{=>} "f";"g" & & = &  & T^2 \ar@{-}[dr]^{\mu}   \ar@{-}[dl]_{1}  \ar@{}[ur]_{1}^(.25){}="s"^(.75){}="r"   \ar@{=>} "s";"r" & &  \\ 
T^2  \ar@{-}[dr]_{\mu}  \ar@{}[ur]^{a}^(.25){}="c"^(.75){}="d"  \ar@{=>} "c";"d" && T \ar@{-}[dl]^{1}   &&  T^2  \ar@{-}[dr]_{\mu} \ar@{}[rr]^{1 }^(.3){}="j"^(.7){}="k"  \ar@{=>} "j";"k"&& T  \ar@{-}[dl]^{1}  \\ & T &&&&T && \\ 
} \]

It is straightforward to verify the following. 

\begin{theorem}
Each two-level adjunction $(F, U, \eta, \epsilon, f u) $ determines a lax-Gray monad $$ (T, \eta, \mu,  a, p, q) = (UF, \eta, U\epsilon F, U \epsilon \epsilon F, Uf, uF)$$
\end{theorem}

\section{Universal Monads, Monoids and the Simplicial Category}

This section presents some familiar results in a form suitable for generalization to lax-Gray categories in the next section.  

Let {\sf 1-Mnd} be the 2-category with 1 object $*$, freely generated by one 1-cell $T$ and two 2-cells $\eta:  1 \to T$ and $\mu:  T^2 \to T$ subject to the usual coherence conditions $\mu. T \mu = \mu . \mu T$ and $\mu . T\eta = 1 = \mu . \eta T$.  Then $(T, \eta, \mu)$ satisfies the monad equations on $*$.  Moreover, $(T, \eta, \mu)$ is a {\bf universal} monad in the sense that for each monad $(T_X, \mu_X, \eta_X)$ on an object $X$ of a 2-category $\mathbb{C}$, there exists a unique functor $F:  {\sf 1-Mnd} \to \mathbb{C}$ with $F(*) = X$ and $(FT, F\eta, F\mu) = (T_X, \eta_X, \mu_X)$.  

Schanual-Street \cite{10} give a concrete description of {\sf 1-Mnd}, which they call {\sf Mnd}.   Their category ${\sf Mnd}$ has one object $0$ and arrows $\sf Mnd(0, 0)$ given by  the category $\Delta$ of finite ordinals and order preserving functions.  Thus as a 2-category, {\sf Mnd} has the finite ordinals as 1-cells and ordinal addition as 1-cell composition.

The ordinal $1$ is terminal in {\sf Mnd} and there are unique 2-cells $\mu:  2 \to 1$ and $\eta:  0 \to 1$, and $(1, \eta, \mu)$ is a monad on $0$ (Kelly-Street \cite{6}) in the 2-category {\sf Mnd}.   Desuspending, we get the classical result from  
 (MacLane \cite{9}) that the 3-tuple $(1, \eta, \mu)$  is a monoid in $\Delta = Mnd(0,0)$ and is universal. 
 
 Now {\sf 1-Mnd} described above is a freely generated 2-category also containing a universal monad $(T, \eta, \mu)$ which under desuspension becomes the universal monoid $(T, \eta, \mu)$ in the symmetric monoidal category $({\sf 1-Mnd})(*, *), \otimes, 1)$.   This may be regarded as an alternate (or abstract) description of $(1, \eta, \mu)$ as a monad in $\Delta$.

This provides us with a model for extending results to the lax-Gray case (MacDonald-Stone \cite{8}). 

\section{The Universal lax-Gray 2-Monad} \label{S:2mnd}

Let {\sf 2-Mnd} denote the lax-Gray category with one object $*$ and freely generated by a single 1-cell $T$, two lax natural 2-cells $\mu:  T^2 \to T$ and $\eta:  1 \to T$ and three 3-cells $a:  \mu T \mu \Rightarrow \mu . \mu T, p: \mu . T \eta \Rightarrow 1_T$ and $q:  1_T \Rightarrow \mu . \eta T$ subject to the coherence conditions [M$a$], [M$\eta$], [M$\mu$], [M$q$] and  [M$p$] drawn in Section 3.  

\begin{lemma}
The 6-tuple $(T, \eta, \mu, a, p, q)$ is a lax-Gray 2-monad on the object $*$ of the lax-Gray category {\sf 2-Mnd}. 

\end{lemma}

\begin{theorem}
The 2-monad $(T, \eta, \mu, a, p, q)$ on $*$ in {\sf 2-Mnd} is universal for lax-Gray 2-monads on objects of a lax-Gray category.
\end{theorem}

In principal, the preceding description of {\sf 2-Mnd} given by generators and coherence conditions completely defines this category.  However, we want to explicitly describe the 1-, 2-, and 3-cells appearing in this freely generated category.

{\sf 2-Mnd} has one object $*$, and the 1-cells are freely generated by $T$, hence explicitly described by $1, T, T^2, T^3, \dots, T^n, \dots$.    The 2-cells are freely generated by $\mu:  T^2 \to T$ and $\eta:  1 \to T$.  Describing these 2-cells explicitly is slightly more complicated.  We start by defining the {\bf basic} 2-cells of {\sf 2-Mnd}, which are created by inserting $\eta$ and $\mu$ into a finite string of $T$'s.  Setting notation, let us denote a cell given by $T^i \eta T^{n-i}$ by $\eta^n_i$, and similarly for $\mu$:  
\[ \xymatrix{  T^n  \ar[dd]^{ \eta^n_i \, \, \,  }  && T^i1T^{n-i}\ar[dd]_{T^i\eta T^{n-i} }  &&& T^{n+1} \ar[dd]^{\mu^n_i }  &&  T^iT^2T^{n-i-1} \ar[dd]_{ \, \, T^i\mu T^{n-i-1}} \\
&=& &&&&=
\\ T^{n+1} &&T^i T T^{n-i}   &&&  T^n &&  T^iTT^{n-i-1}
 \\ }\] 
Then the general 2-cells of {\sf 2-Mnd} are finite strings of basic 2-cells with appropriately matching domains and codomains and obvious composition of strings.  

The 3-cells of {\sf 2-Mnd} are freely generated by $a, p$ and $q$ subject to our five coherence conditions.    These include the generating 3-cells 
\[ \xymatrix{ & T^3 \ar[dl]_{\mu T = \mu_0^2}  \ar[dr]^{\mu_1^2 = T\mu} \\
T^2  \ar[dr]_{\mu_0^1 = \mu }  \ar@{}[rr]^{a }^(.35){}="a"^(.65){}="b"  \ar@{=>} "a";"b"  && T^2 \ar[dl]^{\mu = \mu^1_0}  \\
& T \\  
} \phantom{wwww}    \xymatrix{ & T \ar[dl]_{T \eta }  \ar[dr]^{1} \\
T^2  \ar[dr]_{\mu_0^1 = \mu }  \ar@{}[rr]^{p }^(.35){}="a"^(.65){}="b"  \ar@{=>} "a";"b"  && T \ar[dl]^{1}  \\
& T \\  
}\] 
\[ \xymatrix{ & T \ar[dl]_{1}  \ar[dr]^{\eta T = \eta^1_0} \\
T  \ar[dr]_{1}  \ar@{}[rr]^{q }^(.35){}="a"^(.65){}="b"  \ar@{=>} "a";"b"  && T^2 \ar[dl]^{\mu = \mu^1_0}  \\
& T \\  
}\] 

In addition, we have 3-cells $\eta \eta, \mu \mu, \mu \eta$ and $\eta \mu$ obtained from the lax 2-cells $\eta$ and $\mu$.    To see how these cells get generated, we first look at the simpler case of {\sf Cat}.  If $\eta $ is a lax-natural transformation, then for each 1-cell $\alpha:  X \to Y$ we have a 2-cell 
 \[ \xymatrix{ X\ar[r]^{\eta_X} \ar[d]_{\alpha} & TX \ar[d]^{T\alpha}  \\ 
     Y \ar[r]_{\eta_Y} \ar@{}[ur]^{\eta_\alpha}^(.3){}="a"^(.7){}="b"  \ar@{=>} "a";"b"  & TY } \] 
and thus for each object $X$, there is a 2-cell
 \[ \xymatrix{ X\ar[r]^{\eta_X} \ar[d]_{\eta_X} & TX \ar[d]^{T\eta_X}  \\ 
     TX \ar[r]_{\eta_{TX}} \ar@{}[ur]^{\eta_{\eta_X}}^(.3){}="a"^(.7){}="b"  \ar@{=>} "a";"b"  & T^2X } \] 
and these 2-cells determine a modification 
 \[ \xymatrix{ I\ar[r]^{\eta} \ar[d]_{\eta} & T \ar[d]^{T\eta}  \\ 
     T \ar[r]_{\eta_{T}} \ar@{}[ur]^{\eta{\eta}}^(.3){}="a"^(.7){}="b"  \ar@{=>} "a";"b"  & T^2 }  \phantom{ww} \textup{ or } \phantom{ww} 
     \xymatrix{ & T \ar[dl]_{\eta_0^0}  \ar[dr]^{\eta_0^0} \\
T  \ar[dr]_{T \eta = \eta_1^1}  \ar@{}[rr]^{\eta \eta  }^(.35){}="a"^(.65){}="b"  \ar@{=>} "a";"b"  && T \ar[dl]^{\eta^1_0 = \eta T}  \\
& T \\  
}\] 
where the second diagram presents the modification between our basic 2-cells.  

Similarly for each $X$ there is a 2-cell 
 \[ \xymatrix{ T^2X\ar[r]^{\eta_{T^2X}} \ar[d]_{\mu_X} & T^3X \ar[d]^{T\mu_X}  \\ 
     TX \ar[r]_{\eta_{TX}} \ar@{}[ur]^{\eta_{\mu_X}}^(.3){}="a"^(.7){}="b"  \ar@{=>} "a";"b"  & T^2X } \] 
     which determines the modification 
\[   \xymatrix{ & T^2 \ar[dl]_{\mu_0^1 = \mu}  \ar[dr]^{\eta T^2 = \eta^2_0} \\
T  \ar[dr]_{T \eta^1_0 = \eta T}  \ar@{}[rr]^{\eta \mu }^(.35){}="a"^(.65){}="b"  \ar@{=>} "a";"b"  && T^3 \ar[dl]^{T \mu = \mu^2_1}  \\
& T^2 \\  
}\] 

In the same way, the lax-natural transformation $\mu$ gives, for each 1-cell $\alpha:  X \to Y$, a 2-cell 
\[ \xymatrix{ T^2X\ar[d]_{T^2\alpha} \ar[r]^{\mu_X} & TX \ar[d]^{T\alpha}  \\ 
     T^2Y \ar[r]_{\mu_{Y}} \ar@{}[ur]^{\mu_{\alpha}}^(.3){}="a"^(.7){}="b"  \ar@{=>} "a";"b"  & T Y } \] 
     Thus for each object $X$ we have 2-cells 
     \[ \xymatrix{ T^4X\ar[d]_{T^2\mu_X} \ar[r]^{\mu_{T^2X}} & T^3X \ar[d]^{T\mu_X}  \\ 
     T^3X \ar[r]_{\mu_{TX}} \ar@{}[ur]^{\mu_{\mu_X}}^(.3){}="a"^(.7){}="b"  \ar@{=>} "a";"b"  & T^2X } \phantom{ww} \textup{ and } \phantom{ww}\xymatrix{ T^2X\ar[d]_{T^2\eta} \ar[r]^{\mu_X} & TX \ar[d]^{T\eta}  \\ 
     T^3X \ar[r]_{\mu_{TX}} \ar@{}[ur]^{\mu_{\eta_X}}^(.3){}="a"^(.7){}="b"  \ar@{=>} "a";"b"  & T^2X } \] with corresponding modifications 
     \[   \xymatrix{ & T^4 \ar[dl]_{T^2 \mu}  \ar[dr]^{\mu_{T^2} } \\
T^3  \ar[dr]_{\mu_{T}}  \ar@{}[rr]^{\mu \mu }^(.35){}="a"^(.65){}="b"  \ar@{=>} "a";"b"  && T^3 \ar[dl]^{T \mu }  \\
& T^2 \\  
}\phantom{ww}  \textup{ and } \phantom{ww} \xymatrix{ & T^2 \ar[dl]_{T^2 \eta}  \ar[dr]^{\mu} \\
T^3  \ar[dr]_{\mu_T}  \ar@{}[rr]^{\mu \eta  }^(.35){}="a"^(.65){}="b"  \ar@{=>} "a";"b"  && T \ar[dl]^{T \eta}  \\
& T^2 \\  
}\] 

Thus we can describe generating 3-cells as follows:  for each generating 3-cell $\zeta$, we can define $\zeta^n_i = T^i \zeta T^{n-1}$.    So for example we have cells $a_i^n = T^i a T^{n-1}$ as pictured:  
\[ \xymatrix{ & T^i T^3 T^{n-i} \ar[dl]_{\mu_{i+1}^{n+2} = T^i T\mu T^{n-i}}  \ar[dr]^{T^i \mu T T^{n-i}= \mu_i^{n+2}} \\
T^i T^2 T^{n-i}  \ar[dr]_{\mu_i^{n+1} = T^i \mu T^{n-i}}  \ar@{}[rr]^{a^n_i }^(.35){}="a"^(.65){}="b"  \ar@{=>} "a";"b"  && T^i T^2 T^{n-i} \ar[dl]^{T^i \mu T^{n-i} = \mu_i^{n+1}}  \\
& T^i T T^{n-i}   
}\] 
and similarly $p^n_i = T^i p T^{n-i}, q^n_i = T^i q T^{n-i}, \eta \eta_n^i = T^i \eta \eta T^{n-1}$ etc.  
Then the 3-cells are finite strings of basic 3-cells with matching domains and codomains.

\section{The Lax-Gray Simplicial Category $\Delta_{LG}$}

We now define the higher dimensional analogue of the simplicial category $\Delta$.   
The structure of this category $\Delta_{LG}$ will mirror the structure of the universal lax-Gray 2-monad, with everything desuspended down one dimension.

Let {\sf 2-Mnd} be the one object lax-Gray category described in Section \ref{S:2mnd}.  So by definition \[ \Delta_{LG} = ({\sf 2-Mnd}(*, *), \otimes, 1_*)\]  Explicitly $\Delta_{LG}$ has objects corresponding to 1-cells $1, T, T^2, \dots$ of {\sf 2-Mnd}, and arrows corresponding to  2-cells given by finite strings of composable basic cells $\eta:  1 \to T, \mu:  T^2 \to T$, the basic cells $\eta^n_i = T^i \eta T^{n-i}$ and $\mu^n_i= T^i \mu T^{n-i}$.    Thus the structure here is the same as what we have in $\Delta$, with objects indexed by positive integers, and arrows generated by face and degeneracy maps. 

However, in $\Delta_{LG}$, these arrows do not satisfy the coherence identities of $\Delta$.  Instead, we have 2-cells connecting the various combinations.  These 2-cells of $\Delta_{LG}$ correspond to  the 3-cells of {\sf 2-Mnd}.   For example, in {\sf 1-Mnd}, there is an associativity identity that says that $\mu. T \mu = \mu . \mu T $.  In {\sf 2-Mnd}, this is replaced by an an associativity   3-cell $a:  \mu.T\mu \Rightarrow \mu.\mu T$.   This 3-cell appears one level down in $\Delta_{LG}$ as a 2-cell connecting the arrows  $\mu. T\mu  \Rightarrow \mu.\mu T$. Similarly, the coherence condition on the 3-cell $a$ (an equation involving the  two composites of 3-cells given in the previous section) becomes in $\Delta_{LG}$ a coherence condition with the exact same equation, but now involving 2-cells.  

The 3-cells of  {\sf 2-Mnd} were analyzed in the previous section, and are defined from generating cells $a, p, q$ and the cells $\eta \eta, \mu \mu, \mu \eta$ and $\eta \mu$ obtained from the lax 2-cells $\eta$ and $\mu$.  From these, we can create finite strings of basic 2-cells for $\Delta_{LG} $, of the form $a^n_i = T^i a T^{n-i}$ and similarly built on  $p, q, \eta \eta, \mu \mu, \mu \eta$ and $\eta \mu$ as described in the previous section.


\end{document}